\documentclass[11pt]{article}

\usepackage{latexsym}
\usepackage{amssymb}

\newtheorem{theorem}{Theorem}
\newtheorem{proposition}{Proposition}
\newtheorem{corollary}{Corollary}
\newtheorem{definition}{Definition}
\newtheorem{lemma}{Lemma}

\newtheorem{example}{Example}

\begin{document}
{
\begin{center}
{\Large\bf
Difference equations related to Jacobi-type pencils.}
\end{center}
\begin{center}
{\bf S.M. Zagorodnyuk}
\end{center}
\section{Introduction.}

%%%%%%%%%%%%%%%%%%%%%%%%%%%%%%%%%%   17:35  9.02.2018

The theory of orthogonal polynomials on the real line has a lot of old and new contributions and applications, see~\cite{cit_50000_Gabor_Szego},
\cite{cit_20000_Suetin}, \cite{cit_3000_Chihara},~\cite{cit_5000_Ismail}. 
This theory is closely related to the theory of Jacobi matrices, by using of the associated second order linear difference equation as a link.
The classical theory of orthogonal polynomials has various generalizations, e.g., we can mention the matrix and operator orthogonal
polynomials.
In this paper we shall deal with orthogonal polynomials associated with Jacobi-type pencils~\cite{cit_95000_Z},\cite{cit_97000_Z}.
In particular, these orthogonal polynomials are related to a discrete grid model for a certain linear pencil of fourth-order
differential operators. The latter linear pencil of differential operators appears in some physical applications and it was recently studied by
Ben Amara, Vladimirov and Shkalikov, see~\cite{cit_8_Ben_Amara_Vladimirov_Shkalikov}.

Observe that linear pencils of Jacobi matrices lead to (bi)orthogonal rational functions. Such pencils and their properties are
intensively studied nowadays, see~\cite{cit_98000_Zhedanov_JAT},\cite{cit_98100_Dereviagin} and references therein.  

Recall the following definition from~\cite{cit_95000_Z}.

\begin{definition}
\label{d1_1}
A set $\Theta = \left(
J_3, J_5, \alpha, \beta
\right)$,
where $\alpha>0$, $\beta\in\mathbb{R}$, $J_3$ is a Jacobi matrix and
$J_5$ is a semi-infinite real symmetric five-diagonal matrix with positive numbers on the second subdiagonal,
is said to be
\textbf{a Jacobi-type pencil (of matrices)}.
\end{definition}
From this definition we see that matrices $J_3$ and $J_5$ have the following form:
\begin{equation}
\label{f1_5}
J_3 =
\left(
\begin{array}{cccccc}
b_0 & a_0 & 0 & 0 & 0 & \cdots\\
a_0 & b_1 & a_1 & 0 & 0 & \cdots\\
0 & a_1 & b_2 & a_2 & 0 &\cdots\\ 
\vdots & \vdots & \vdots & \ddots \end{array}
\right),\qquad a_k>0,\ b_k\in\mathbb{R},\ k\in\mathbb{Z}_+;
\end{equation}

\begin{equation}
\label{f1_10}
J_5 =
\left(
\begin{array}{ccccccc}
\alpha_0 & \beta_0 & \gamma_0 & 0 & 0 & 0 & \cdots\\
\beta_0 & \alpha_1 & \beta_1 & \gamma_1 & 0 & 0 & \cdots\\
\gamma_0 & \beta_1 & \alpha_2 & \beta_2 & \gamma_2 & 0 & \cdots\\ 
0 & \gamma_1 & \beta_2 & \alpha_3 & \beta_3 & \gamma_3 & \cdots \\
\vdots & \vdots & \vdots &\vdots & \ddots \end{array}
\right),\ \alpha_n,\beta_n\in\mathbb{R},\ \gamma_n>0,\ n\in\mathbb{Z}_+.
\end{equation}

With a Jacobi-type pencil of matrices $\Theta$ one associates a system of polynomials
$\{ p_n(\lambda) \}_{n=0}^\infty$, such that
\begin{equation}
\label{f1_15}
p_0(\lambda) = 1,\quad p_1(\lambda) = \alpha\lambda + \beta,
\end{equation}
and
\begin{equation}
\label{f1_20}
(J_5 - \lambda J_3) \vec p(\lambda) = 0,
\end{equation}
where $\vec p(\lambda) = (p_0(\lambda), p_1(\lambda), p_2(\lambda),\cdots)^T$. Here the superscript $T$ means the transposition.
Polynomials $\{ p_n(\lambda) \}_{n=0}^\infty$ are said to be \textbf{associated to the Jacobi-type pencil of matrices $\Theta$}.

In particular, for each system of orthonormal polynomials on the real line (with positive leading
coefficients and $p_0=1$) one may choose $J_3$ to be the corresponding Jacobi matrix
(which elements are the recurrence coefficients),  $J_5 = J_3^2$, and $\alpha,\beta$  being the coefficients of $p_1$
($p_1(\lambda) = \alpha\lambda + \beta$).

One can rewrite relation~(\ref{f1_20}) in the scalar form:
$$ \gamma_{n-2} p_{n-2}(\lambda) + (\beta_{n-1}-\lambda a_{n-1}) p_{n-1}(\lambda) + (\alpha_n-\lambda b_n) p_n(\lambda) +
$$
\begin{equation}
\label{f1_30}
+ (\beta_n-\lambda a_n) p_{n+1}(\lambda) + \gamma_n p_{n+2}(\lambda) = 0,\qquad n\in\mathbb{Z}_+,
\end{equation}
where $p_{-2}(\lambda) = p_{-1}(\lambda) = 0$, $\gamma_{-2} = \gamma_{-1} = a_{-1} = \beta_{-1} = 0$.

In this paper we shall study various difference equations related to Jacobi-type pencils. Let us briefy describe the content of 
the present paper.
In the next section, for the convenience of the reader, we recall some basic definitions and results from~\cite{cit_95000_Z},
\cite{cit_97000_Z} which will be used in what follows.
In Section~3 we investigate a 4-th order linear difference equation generated by the recurrence relation~(\ref{f1_30}).
A basic set of solution is constructed. An analog of the Christoffel-Darboux formula is given. Some spectral properties
of the corresponding truncated pencil are investigated.

Section~4 is devoted to a special perturbation of orthonormal polynomials on the real line. In this case, the associated
polynomials $p_n(\lambda)$ have a transparent explicit representation. Some special matrix orthonormality relations hold for
$p_n(\lambda)$. These polynomials satisfy a non-standard three term recurrence relation.
Differential equations for classical orthogonal polynomials on the real line imply some differential
equations for the corresponding perturbed polynomials $p_n(\lambda)$. It turnes out that the perturbation $y_n(t)$ of Jacobi polynomials satisfies
the following $4$-th order differential equation with polynomial coefficients (not depending on $n$):
$$ c(t) y^{(4)}(t) + d(t) y'''(t) + f(t) y''(t) + g(t) y'(t) + h(t) y(t) + $$
\begin{equation}
\label{f1_150}
+ \lambda \left(
\varphi(t) y''(t) + \psi(t) y'(t) + \theta(t) y(t)
\right) = 0,
\end{equation}
where $\lambda=\lambda_n\in\mathbb{R}$.
The polynomial $y_n$ is a unique, up to a constant multiple, $n$-th degree real polynomial solution of equation~(\ref{f1_150}).
Moreover, the constant $\lambda$ is uniquely determined. Thus, we have a transparent analogy with classical systems of orthogonal 
polynomials on the real line. 
Notice that orthogonal polynomials, satisfying $4$-th order differential equations like in~(\ref{f1_150}), but with
$\varphi(t) = \psi(t) = 0$, appeared in the literature, see the book~\cite{cit_5100_Krall} and historical references therein.
On the other hand, polynomial solutions for second order differential operators, like in~(\ref{f1_150}) but with
$c(t) = d(t) = 0$, were studied in~\cite{cit_15000_Sawyer},\cite{cit_1000_Chaundy}.

Section~5 is devoted to the following question: \textit{when the associated polynomials $p_n(\lambda)$ are related to
a single block Jacobi matrix and matrix orthogonal polynomials}? It is known that a block Jacobi matrix generates
matrix orthonormal polynomials as well as scalar orthonormal polynomials on radial rays, see, e.g.~\cite{cit_98300_D_A},~\cite{cit_90000_Z}. 
The latter systems of polynomials
have a simple algebraic interrelation. However, the scalar polynomials have different properties of zeros~\cite{cit_94000_Z}.
We shall show that the associated polynomials $p_n(\lambda)$ are related to
a single block Jacobi matrix if and only if the associated operator $\mathcal{A}$ of the pencil has a symmetric power.
For example, we show that some polynomials from Section~4 are not related to a single block Jacobi matrix.

\noindent
{\bf Notations. }
As usual, we denote by $\mathbb{R}, \mathbb{C}, \mathbb{N}, \mathbb{Z}, \mathbb{Z}_+$,
the sets of real numbers, complex numbers, positive integers, integers and non-negative integers,
respectively. By $\mathbb{P}$ we denote the set of all polynomials with complex coefficients.

\noindent
By $l_2$ we denote the usual Hilbert space of all complex sequences $c = (c_n)_{n=0}^\infty = (c_0,c_1,c_2,...)^T$ with the finite norm
$\| c \|_{l_2} = \sqrt{\sum_{n=0}^\infty |c_n|^2}$. Here $T$ means the transposition.
The scalar product of two sequences $c = (c_n)_{n=0}^\infty, d = (d_n)_{n=0}^\infty\in l_2$ is given by
$(c,d)_{l_2} = \sum_{n=0}^\infty c_n \overline{d_n}$. 
We denote $\vec e_m = (\delta_{n,m})_{n=0}^\infty\in l_2$, $m\in\mathbb{Z}_+$.
By $l_{2,fin}$ we denote the set of all finite vectors from $l_2$, i.e. vectors with all, but finite number, components being zeros.

\noindent
By $\mathfrak{B}(\mathbb{R})$ we denote the set of all Borel subsets of $\mathbb{R}$.
If $\sigma$ is a (non-negative) bounded measure on $\mathfrak{B}(\mathbb{R})$ then by $L^2_\sigma$ we denote a Hilbert space
of all (classes of equivalences of) complex-valued functions $f$ on $\mathbb{R}$ with a finite norm
$\| f \|_{L^2_\sigma} = \sqrt{ \int_\mathbb{R} |f(x)|^2 d\sigma }$.
The scalar product of two functions $f,g\in L^2_\sigma$ is given by
$(f,g)_{L^2_\sigma} = \int_{\mathbb{R}} f(x) \overline{g(x)} d\sigma$. 
By $[ f ] = [f]_{L^2_\sigma}$  we denote the class of equivalence in $L^2_\sigma$ which contains the representative $f$.
By $\mathcal{P}$ we denote a set of all (classes of equivalence which contain) complex polynomials in $L^2_\sigma$.
As usual, we sometimes use the representatives instead of their classes in formulas. 
Let $B$ be an arbitrary linear operator in $L^2_\sigma$ with the domain $\mathcal{P}$.
Let $f(\lambda)\in\mathbb{P}$ be nonzero and of degree $d\in\mathbb{Z}_+$,
$f(\lambda) = \sum_{k=0}^d d_k \lambda^k$, $d_k\in\mathbb{C}$. We set
$$ f(B) = \sum_{k=0}^d d_k B^k;\quad
B^0 := E|_{\mathcal{P}}. $$
If $f\equiv 0$, then $f(B):= 0|_{\mathcal{P}}$.

\noindent
If H is a Hilbert space then $(\cdot,\cdot)_H$ and $\| \cdot \|_H$ mean
the scalar product and the norm in $H$, respectively.
Indices may be omitted in obvious cases.
For a linear operator $A$ in $H$, we denote by $D(A)$
its  domain, by $R(A)$ its range, and $A^*$ means the adjoint operator
if it exists. If $A$ is invertible then $A^{-1}$ means its
inverse. $\overline{A}$ means the closure of the operator, if the
operator is closable. If $A$ is bounded then $\| A \|$ denotes its
norm.
For a set $M\subseteq H$
we denote by $\overline{M}$ the closure of $M$ in the norm of $H$.
By $\mathop{\rm Lin}\nolimits M$ we mean
the set of all linear combinations of elements of $M$,
and $\mathop{\rm \overline{span}}\nolimits M := \overline{ \mathop{\rm Lin}\nolimits M }$.
By $E=E_H$ ($0=0_H$) we denote the identity operator in $H$, i.e. $E_H x = x$,
$x\in H$ (respectively the null operator in $H$, i.e. $0_H x = 0$,
$x\in H$). If $H_1$ is a subspace of $H$, then $P_{H_1} =
P_{H_1}^{H}$ is an operator of the orthogonal projection on $H_1$
in $H$.

\section{Preliminaries.}

In this section we recall basic definitions and results from~\cite{cit_95000_Z}, \cite{cit_97000_Z} which will be used later. 
Set
\begin{equation}
\label{f3_40}
u_n := J_3 \vec e_n = a_{n-1} \vec e_{n-1} + b_n \vec e_n + a_n \vec e_{n+1},
\end{equation}
\begin{equation}
\label{f3_50}
w_n := J_5 \vec e_n = \gamma_{n-2} \vec e_{n-2} + \beta_{n-1} \vec e_{n-1} + \alpha_n \vec e_n + \beta_n \vec e_{n+1}
+ \gamma_n \vec e_{n+2},\qquad n\in\mathbb{Z}_+. 
\end{equation}
Here and in what follows by $\vec e_k$ with negative $k$ we mean (vector) zero.
The following operator:   
$$ A f = \frac{\zeta}{\alpha} (\vec e_1 - \beta \vec e_0)
+ 
\sum_{n=0}^\infty \xi_n w_n, $$
\begin{equation}
\label{f3_60}
f = \zeta \vec e_0 + \sum_{n=0}^\infty \xi_n u_n\in l_{2,fin},\quad \zeta, \xi_n\in\mathbb{C}, 
\end{equation}
with $D(A) = l_{2,fin}$
is said to be \textbf{the associated operator for the Jacobi-type pencil $\Theta$}.
Notice that in the sums in~(\ref{f3_60}) only finite number of $\xi_n$ are nonzero. We shall always assume this in the case
of elements from the linear span.
In particular, we have
$$ A J_3 \vec e_n = J_5 \vec e_n,\qquad n\in\mathbb{Z}_+, $$
and therefore
\begin{equation}
\label{f3_60_1}
A J_3 = J_5. 
\end{equation}
As usual, the matrices $J_3$ and $J_5$ define linear operators with the domain $l_{2,fin}$ which we denote by the same letters.

For an arbitrary non-zero polynomial $f(\lambda)\in\mathbb{P}$ of degree $d\in\mathbb{Z}_+$,
$f(\lambda) = \sum_{k=0}^d d_k \lambda^k$, $d_k\in\mathbb{C}$, we set
$f(A) = \sum_{k=0}^d d_k A^k$.
Here $A^0 := E|_{l_{2,fin}}$.
For $f(\lambda) \equiv 0$, we set
$f(A) = 0|_{l_{2,fin}}$.
%The following relations hold~\cite{cit_95000_Z}:
%\begin{equation}
%\label{f3_110}
%\vec e_n = p_n(A) \vec e_0,\qquad n\in\mathbb{Z}_+;
%\end{equation}
%\begin{equation}
%\label{f3_120}
%\left( p_n(A) \vec e_0, p_m(A) \vec e_0 \right)_{l_2} = \delta_{n,m},\qquad n,m\in\mathbb{Z}_+.
%\end{equation}
Denote by $\{ r_n(\lambda) \}_{n=0}^\infty$, $r_0(\lambda) = 1$, the system of polynomials satisfying
\begin{equation}
\label{f3_130}
J_3 \vec r(\lambda) = \lambda \vec r(\lambda),\quad \vec r(\lambda) = (r_0(\lambda), r_1(\lambda),r_2(\lambda), ...)^T.
\end{equation}
These polynomials are orthonormal on the real line with respect to a (possibly non-unique) non-negative finite measure $\sigma$ on $\mathfrak{B}(\mathbb{R})$
(Favard's theorem). 
Consider the following operator:
\begin{equation}
\label{f3_140}
U \sum_{n=0}^\infty \xi_n \vec e_n = \left[ \sum_{n=0}^\infty \xi_n r_n(x) \right ],\qquad \xi_n\in\mathbb{C},
\end{equation}
which maps $l_{2,fin}$ onto $\mathcal{P}$. Here by $\mathcal{P}$ we denote a set of all
(classes of equivalence which contain) complex polynomials in $L^2_\sigma$.
Denote
\begin{equation}
\label{f3_150}
\mathcal{A} = \mathcal{A}_\sigma = U A U^{-1}. 
\end{equation}
The operator $\mathcal{A} = \mathcal{A}_\sigma$ is said to be \textbf{the model representation in $L^2_\sigma$ of the associated operator $A$}.

\begin{theorem} (\cite{cit_95000_Z})
\label{t3_1}
Let $\Theta = \left(
J_3, J_5, \alpha, \beta
\right)$ be a Jacobi-type pencil.
Let $\{ r_n(\lambda) \}_{n=0}^\infty$, $r_0(\lambda) = 1$, be a system of polynomials satisfying~(\ref{f3_130}) and
$\sigma$ be  their (arbitrary) orthogonality measure on $\mathfrak{B}(\mathbb{R})$.
The associated polynomials $\{ p_n(\lambda) \}_{n=0}^\infty$ satisfy the following relations:
\begin{equation}
\label{f3_160}
\int_{\mathbb{R}} p_n(\mathcal{A})(1) \overline{ p_m(\mathcal{A})(1) } d\sigma
= \delta_{n,m},\qquad n,m\in\mathbb{Z}_+,
\end{equation}
where $\mathcal{A}$ is the model representation  in $L^2_\sigma$  of the associated operator $A$.
\end{theorem}

\begin{definition} (\cite{cit_97000_Z})
\label{dd1_5}
Let $\Theta = \left(
J_3, J_5, \alpha, \beta
\right)$ be a Jacobi-type pencil and 
$\{ p_n(\lambda) \}_{n=0}^\infty$ be the associated polynomials to $\Theta$.
A sesquilinear functional $S(u,v)$, $u,v\in\mathbb{P}$, satisfying the following relation:
\begin{equation}
\label{ff1_6}
S(p_n,p_m) = \delta_{n,m},\qquad n,m\in\mathbb{Z}_+,
\end{equation}
is said to be \textbf{the spectral function of the Jacobi-type pencil $\Theta$}.
\end{definition}

The spectral function of the Jacobi-type pencil satisfy the following properties:
\begin{equation}
\label{ff1_7_0}
\overline{S(u,v)} = S(v,u),\qquad u,v\in\mathbb{P};
\end{equation}
and
\begin{equation}
\label{ff1_7_0_1}
S(u,u)\geq 0,\qquad u\in\mathbb{P}.
\end{equation}

\begin{theorem} (\cite{cit_97000_Z})
\label{tt2_1}
A sesquilinear
functional $S(u,v)$, $u,v\in\mathbb{P}$, 
satisfying relations~(\ref{ff1_7_0}), (\ref{ff1_7_0_1}),
is
the spectral function of a Jacobi-type pencil if and only if 
it admits the following integral representation:
\begin{equation}
\label{ff2_5}
S(u,v) = \int_{\mathbb{R}} u(\mathcal{A})(1) \overline{ v(\mathcal{A})(1) } d\sigma,\qquad u,v\in\mathbb{P},
\end{equation}
where $\sigma$ is a non-negative measure on $\mathfrak{B}(\mathbb{R})$ with all finite power moments, 
$\int_\mathbb{R} d\sigma = 1$, $\int_\mathbb{R} |g(x)|^2 d\sigma > 0$, for any non-zero complex polynomial $g$, 
and $\mathcal{A}$ is a linear operator in $L^2_\sigma$ with the following properties:
\begin{itemize}
\item[(i)] $D(\mathcal{A}) = \mathcal{P}$;

\item[(ii)] For each $k\in\mathbb{Z}_+$ it holds:
\begin{equation}
\label{ff2_7}
\mathcal{A} x^k = \xi_{k,k+1} x^{k+1} + \sum_{j=0}^k \xi_{k,j} x^j,
\end{equation}
where $\xi_{k,k+1} > 0$, $\xi_{k,j}\in\mathbb{R}$ ($0\leq j\leq k$);

\item[(iii)] The operator $\mathcal{A} \Lambda_0$ is symmetric. 
Here by $\Lambda_0$ we denote the operator of the multiplication by an independent variable in $L^2_\sigma$
restricted to $\mathcal{P}$.
\end{itemize}
\end{theorem}

\begin{corollary} (\cite{cit_97000_Z})
\label{cc2_1}
Let $\sigma$ be a non-negative measure on $\mathfrak{B}(\mathbb{R})$ with all finite power moments, 
$\int_\mathbb{R} d\sigma = 1$, $\int_\mathbb{R} |g(x)|^2 d\sigma > 0$, for any non-zero complex polynomial $g$. 
A linear operator $\mathcal{A}$ in $L^2_\sigma$ is a model representation in $L^2_\sigma$ of the 
associated operator of a Jacobi-type pencil if and only if properties $(i)$-$(iii)$ of Theorem~\ref{tt2_1} hold.
\end{corollary}

\section{The fourth order linear difference equation related to a Jacobi-type pencil.}

Let $\Theta = \left( J_3, J_5, \alpha, \beta \right)$ be a Jacobi type pencil with matrices $J_3, J_5$ having form~(\ref{f1_5}),(\ref{f1_10}).
Consider the following $4$-th order linear difference equation with complex unknowns $\{ y_k \}_{k=0}^\infty$:
$$ \gamma_{n-2} y_{n-2} + (\beta_{n-1}-\lambda a_{n-1}) y_{n-1} + (\alpha_n-\lambda b_n) y_n +
(\beta_n-\lambda a_n) y_{n+1} +
$$
\begin{equation}
\label{f3_5}
+ \gamma_n y_{n+2} = 0,\qquad n = 2,3,....
\end{equation}

Let us construct the basic set of solutions for equation~(\ref{f3_5}).
Of course, the associated (to $\Theta$) polynomials $\{ p_n(\lambda) \}_{n=0}^\infty$ give a solution to equation~(\ref{f3_5}). 
Let $S(u,v)$ be the spectral function for the pencil $\Theta$. Then we may introduce the following system of polynomials:
\begin{equation}
\label{f3_7}
q_n(\lambda) = S_t \left(
\frac{p_n(\lambda) - p_n(t)}{\lambda - t}, 1
\right),\qquad n\in\mathbb{Z}_+,
\end{equation}
where the subscript $t$ means that the arguments of $S$ are polynomials of $t$.
It is natural to call $q_n(\lambda)$ \textbf{the polynomials of the second kind for the Jacobi-type pencil $\Theta$}.
A direct calculation, with a use of the recurrence relation~(\ref{f1_30}), shows that
$$ \gamma_{n-2} q_{n-2}(\lambda) + (\beta_{n-1}-\lambda a_{n-1}) q_{n-1}(\lambda) + (\alpha_n-\lambda b_n) q_n(\lambda) +
(\beta_n-\lambda a_n) q_{n+1}(\lambda) +
$$
\begin{equation}
\label{f3_9}
+ \gamma_n q_{n+2}(\lambda) = 
S_t\left(a_{n-1}p_{n-1}(t) + b_n p_n(t) + a_n p_{n+1}(t), 1\right)
=
\left\{
\begin{array}{ccc} 0, & n=2,3,...\\
a_0, & n=1\\
b_0, & n=0 \end{array}
\right..
\end{equation}
Thus, $q_n$ give another solution of~(\ref{f3_5}). 
In order to construct more solutions of difference equation~(\ref{f3_5}) we need the following definition.
\begin{definition}
\label{d3_1}
Let $\Theta = \left( J_3, J_5, \alpha, \beta \right)$ be a Jacobi-type pencil, with matrices $J_3, J_5$ as in~(\ref{f1_5}),(\ref{f1_10}).
Let $J_3',J_5'$ be semi-infinite matrices obtained respectively from $J_3,J_5$, by removing of the first row and the first column.
The Jacobi-type pencil $\Theta' = \left( J_3', J_5', \frac{a_0}{\gamma_0}, -\frac{\beta_0}{\gamma_0} \right)$ is said to be \textbf{the shifted pencil
(for the pencil $\Theta$)}.
\end{definition}

Consider the shifted pencil $\Theta' = \left( J_3', J_5', \frac{a_0}{\gamma_0}, -\frac{\beta_0}{\gamma_0} \right)$ for $\Theta$.
Denote by $f_n(\lambda)$ ($n\in\mathbb{Z}_+$) the associated polynomials for $\Theta'$. 
Set
\begin{equation}
\label{f3_10}
u_0(\lambda) = 0,\ u_k(\lambda) = f_{k-1}(\lambda),\ k\in\mathbb{N}.
\end{equation}
By the recurrence relation for $f_n$ it follows that $\{ u_n \}_{n=0}^\infty$ satisfy relation~(\ref{f3_5}) for $n=1,2,...$.
On the other hand, the choice of the constants in $\Theta'$ provides that
$\{ u_n \}_{n=0}^\infty$ satisfy relation~(\ref{f3_5}) for $n=0$.
We shall call $\{ u_n(\lambda) \}_{n=0}^\infty$ \textbf{the shifted polynomials for the pencil $\Theta$}.

\noindent
In a similar manner, consider the shifted pencil for $\Theta'$:
$$ \widetilde\Theta := \left(\Theta'\right)'=
\left( (J_3')',(J_5')', \frac{a_1}{\gamma_1}, -\frac{\beta_1}{\gamma_1} \right). $$ 
Denote by $\widetilde f_n(\lambda)$ ($n\in\mathbb{Z}_+$) the associated polynomials for $\widetilde\Theta$ and 
set
\begin{equation}
\label{f3_15}
w_0(\lambda) = 0,\ w_1(\lambda) = 0,\  w_k(\lambda) = \widetilde f_{k-2}(\lambda),\ k=2,3,....
\end{equation}
By the recurrence relation for $\widetilde f_n$ we conclude that
$\{ w_n \}_{n=0}^\infty$ is a solution of difference equation~(\ref{f3_5}) (for $n\geq 2$).
The choice of the constants in $\Theta''$ provides that
$\{ w_n \}_{n=0}^\infty$ satisfy relation~(\ref{f3_5}) for $n=1$.
The polynomials $\{ w_n(\lambda) \}_{n=0}^\infty$ are said to be \textbf{the double-shifted polynomials for the pencil $\Theta$}.

Thus, using the associated polynomials for some Jacobi-type pencils and the spectral function we obtained four solutions to the difference
equation~(\ref{f3_5}). Below we shall show that these solutions are linearly independent.

We should remark that the construction of the associated orthogonal polynomials is not an easy matter. However, there
exists at least a determinant representation. In fact, consider the moments of the spectral function
$S(u,v)$ of the pencil $\Theta$:
\begin{equation}
\label{f3_17}
s_{m,n} := S(\lambda^m,\lambda^n),\qquad m,n\in\mathbb{Z}_+.
\end{equation}
By the integral representation~(\ref{ff2_5}) and the fact that $\mathcal{A}(1)$ is a real polynomial (see the formula
after~(3.7) in~\cite[p. 9]{cit_97000_Z}) it follows that $s_{m,n}$ are real, and $s_{m,n} = s_{n,m}$.
The spectral function $S(u,v)$ defines an inner product on the complex vector space $\mathbb{P}$, see~\cite[p. 9]{cit_97000_Z}.
In particular, the property
\begin{equation}
\label{f3_19}
\left( S(u,u)=0 \right)\ \Rightarrow\  \left( u=0 \right) 
\end{equation}
holds.
%Denote by $\mathfrak{P}$ the Hilbert space obtained
%by the completion of $\mathbb{P}$ with respect to the corresponding norm.
Choose an arbitrary complex polynomial $u(\lambda)$ of degree $n$:
$$ u(\lambda) = \sum_{j=0}^n \xi_j \lambda^j,\qquad \xi_j\in\mathbb{C}. $$
We may write:
\begin{equation}
\label{f3_22}
0 \leq S(u,u) = \sum_{j,k=0}^n \xi_j \overline{\xi_k} s_{j,k}. 
\end{equation}
By~(\ref{f3_19}) we see that the quadratic form in~(\ref{f3_22}) is positive. Therefore
\begin{equation}
\label{f3_24}
\Delta_n := \left| 
s_{j,k}
\right|_{j,k=0}^n =
\left| 
s_{k,j}
\right|_{j,k=0}^n
> 0,\qquad n\in\mathbb{Z}_+;\quad \Delta_{-1} := 1.
\end{equation}
In a standard way (see, e.g.,~\cite{cit_20000_Suetin}) it can be checked that the following polynomials:
$$ \mathbf{p}_0(\lambda) = 1, $$
\begin{equation}
\label{f3_26}
\mathbf{p}_n(\lambda) = \frac{1}{ \sqrt{\Delta_{n-1} \Delta_n} }
\left|
\begin{array}{ccccc} s_{0,0} & s_{1,0} & \ldots & s_{n,0} \\
s_{0,1} & s_{1,1} & \ldots & s_{n,1} \\
\vdots & \vdots & \ddots & \vdots \\
s_{0,n-1} & s_{1,n-1} & \ldots & s_{n,n-1} \\
1 & \lambda & \ldots & \lambda^n \end{array}
\right|,\qquad n\in\mathbb{N},
\end{equation}
are orthonormal with respect to $S$. 
In fact, expanding the determinant for $\mathbf{p_n}$ ($n\in\mathbb{N}$) by the last row we get
$$ S(\mathbf{p}_n(\lambda), \lambda^k) = 
\left\{
\begin{array}{cc}
0, & k=0,1,...,n-1\\
\sqrt{ \frac{ \Delta_n }{ \Delta_{n-1} } }, & k=n
\end{array}
\right.,
$$
and then it remains to use the properties of $S$.

Since $\mathbf{p}_n$ and $p_n$ ($n\in\mathbb{N}$) have positive leading coefficients, then
$\mathbf{p}_n - \xi p_n$ is a real polynomial of degree at most $n-1$, for a suitable $\xi > 0$.
Then we may write:
$$ S(\mathbf{p}_n - \xi p_n, \mathbf{p}_n - \xi p_n) = 
S(\mathbf{p}_n , \mathbf{p}_n - \xi p_n) - \xi S(p_n, \mathbf{p}_n - \xi p_n) = 0. $$
By~(\ref{f3_19}) we obtain that $\mathbf{p}_n = \xi p_n$.
Since 
$$ 1 = S(\mathbf{p}_n,\mathbf{p}_n) = \xi^2 S(p_n,p_n) = \xi^2, $$
then $\xi = 1$.
Thus, $\mathbf{p}_n = p_n$ ($n\in\mathbb{Z}_+$).

\begin{theorem}
\label{tt3_1}
Let $\Theta = \left( J_3, J_5, \alpha, \beta \right)$ be a Jacobi-type pencil, with matrices $J_3, J_5$ as in~(\ref{f1_5}),(\ref{f1_10}).
The basic set of solutions of the difference equation~(\ref{f3_5}) is given by
the associated polynomials $\{ p_n(\lambda) \}_{n=0}^\infty$, the polynomials of the second kind $\{ q_n(\lambda) \}_{n=0}^\infty$,  
the shifted polynomials $\{ u_n(\lambda) \}_{n=0}^\infty$ and the double-shifted polynomials $\{ w_n(\lambda) \}_{n=0}^\infty$.  
\end{theorem}
\textbf{Proof.} It remains to verify that the above four solutions are linearly independent.
Choose arbitrary complex constants $c_1,c_2,c_3,c_4$ such that
\begin{equation}
\label{f3_28}
c_1 p_n (\lambda) + c_2 q_n(\lambda) + c_3 u_n(\lambda) + c_4 w_n(\lambda) = 0,\qquad n\in\mathbb{Z}_+; (\lambda\in\mathbb{C}). 
\end{equation}
Suppose that $c_2\not= 0$. Then
\begin{equation}
\label{f3_30}
q_n(\lambda) = -\frac{c_1}{c_2} p_n(\lambda) - \frac{c_3}{c_2} u_n(\lambda) - \frac{c_4}{c_2} w_n(\lambda),\qquad n\in\mathbb{Z}_+.
\end{equation}
We obtained a contradiction, since the right-hand side of~(\ref{f3_30}) satisfy relation~(\ref{f3_5}) for $n=1$,
while $q_n(\lambda)$ do not satisfy it.
Therefore $c_2 = 0$.
Suppose that $c_1 \not= 0$.  Then
\begin{equation}
\label{f3_33}
p_n(\lambda) = -\frac{c_3}{c_1} u_n(\lambda) - \frac{c_4}{c_1} w_n(\lambda),\qquad n\in\mathbb{Z}_+.
\end{equation}
Choosing $n=0$ we get a contradiction. Thus, $c_1 = 0$. Suppose that $c_3\not= 0$. Then
\begin{equation}
\label{f3_34}
u_n(\lambda) = - \frac{c_4}{c_3} w_n(\lambda),\qquad n\in\mathbb{Z}_+.
\end{equation}
Since $w_1 = 0$ and $u_1 = 1$ we obtain a contradiction. Therefore $c_3 = 0$.
Since $w_2 = 1$, we conclude that $c_4 = 0$. Therefore the above four solutions are linearly independent.
$\Box$

The following analog of the Christoffel-Darboux formula holds.
\begin{theorem}
\label{tt3_2}
Let $\Theta = \left( J_3, J_5, \alpha, \beta \right)$ be a Jacobi-type pencil, with matrices $J_3, J_5$ as in~(\ref{f1_5}),(\ref{f1_10}).
Let $\{ p_n(\lambda) \}_{n=0}^\infty$ be the associated polynomials to $\Theta$.
For arbitrary $n\in\mathbb{N}$ and $\lambda,y\in\mathbb{C}$: $\lambda\not=y$, the following relation holds:
$$ \sum_{k=0}^n (a_{k-1} p_{k-1}(y) + b_k p_k(y) + a_k p_{k+1}(y)) p_k(\lambda) = $$
$$ = \gamma_{n-1}
\frac{ p_{n+1}(\lambda) p_{n-1}(y) - p_{n+1}(y) p_{n-1}(\lambda)  }{\lambda - y} + \gamma_{n}
\frac{ p_{n+2}(\lambda) p_{n}(y) - p_{n+2}(y) p_{n}(\lambda)  }{\lambda - y} + $$
\begin{equation}
\label{f3_36}
+ ( \beta_n - a_n \lambda ) \frac{ p_{n+1}(\lambda) p_{n}(y) - p_{n+1}(y) p_{n}(\lambda)  }{\lambda - y}.
\end{equation}
\end{theorem}
\textbf{Proof.} 
Multiply the recurrence relation~(\ref{f1_30}) by $p_n(y)$, $y\in\mathbb{C}$:
$$ \gamma_{n-2} p_{n-2}(\lambda) p_n(y) + \beta_{n-1} p_{n-1}(\lambda) p_n(y) + \alpha_n p_n(\lambda) p_n(y) + $$
$$  + \beta_n p_{n+1}(\lambda) p_n(y) + \gamma_n p_{n+2}(\lambda) p_n(y) = $$
\begin{equation}
\label{f3_38}
= \lambda ( a_{n-1} p_{n-1}(\lambda) p_n(y) + b_n p_n(\lambda) p_n(y) + a_n p_{n+1}(\lambda) p_n(y) ),\qquad n\in\mathbb{Z}_+.
\end{equation}
Changing the roles of $\lambda$ and $y$ we get
$$ \gamma_{n-2} p_{n-2}(y) p_n(\lambda) + \beta_{n-1} p_{n-1}(y) p_n(\lambda) + \alpha_n p_n(y) p_n(\lambda) + $$
$$  + \beta_n p_{n+1}(y) p_n(\lambda) + \gamma_n p_{n+2}(y) p_n(\lambda) = $$
\begin{equation}
\label{f3_45}
= y ( a_{n-1} p_{n-1}(y) p_n(\lambda) + b_n p_n(y) p_n(\lambda) + a_n p_{n+1}(y) p_n(\lambda) ),\qquad n\in\mathbb{Z}_+.
\end{equation}
Subtracting the above relations we obtain that
$$ -\Gamma_{n-2} - \Phi_{n-1} + \Phi_n + \Gamma_n =
a_{n-1}(\lambda p_{n-1}(\lambda) p_n(y) - y p_{n-1}(y) p_n(\lambda)) + $$
\begin{equation}
\label{f3_47}
+ (\lambda - y) b_n p_n(\lambda) p_n(y) +
a_n(\lambda p_{n+1}(\lambda) p_n(y) - y p_{n+1}(y) p_n(\lambda)),\qquad n\in\mathbb{Z}_+,
\end{equation}
where
$$ \Gamma_n := \gamma_n ( p_{n+2}(\lambda) p_n(y) - p_{n+2}(y) p_n(\lambda) ), $$
$$ \Phi_n := \beta_n ( p_{n+1}(\lambda) p_n(y) - p_{n+1}(y) p_n(\lambda) ), $$
and $\Gamma_{-2} = \Gamma_{-1} = \Phi_{-1} = 0$.
The right-hand side of relation~(\ref{f3_47}) can be written as
\begin{equation}
\label{f3_49}
-A_{n-1} + A_n + (\lambda - y) (a_{n-1}p_{n-1}(y) + b_n p_n(y) + a_n p_{n+1}(y)) p_n(\lambda), 
\end{equation}
where
$$ A_n := a_n \lambda (p_{n+1}(\lambda) p_n(y) - p_{n+1}(y) p_n(\lambda)),\quad n\in\mathbb{Z}_+, $$
and $A_{-1} = 0$.
Summing in~(\ref{f3_47}) for $n=0$ to $n=r$ and changing indexes we obtain relation~(\ref{f3_36}).
$\Box$

Denote by $J_{3;k}$, $J_{5;k}$ the matrices standing on the intersection of the first $(k+1)$ rows and
the first $(k+1)$ columns of matrices $J_3$, $J_5$, respectively ($k\in\mathbb{Z}_+$).
Set
\begin{equation}
\label{f3_51}
D_j(\lambda) := \det \left(
J_{5;j-1} - \lambda J_{3;j-1}
\right),\qquad j\in\mathbb{N};\ D_0(\lambda):= 1;\ \lambda\in\mathbb{C}.
\end{equation}
The polynomial $D_{j}(\lambda)$ is said to be \textit{the characteristic polynomial} of the self-adjoint linear matrix pencil
$\Phi_j = J_{5;j-1} - \lambda J_{3;j-1}$. Its roots need not to be real values. However, if $J_{3;j-1}>0$ they
should be real, see~\cite[p. 337]{cit_80000_Wilkinson}. 
The roots of $D_{j}(\lambda)$ are said to be \textit{the eigenvalues of the matrix pencil $\Phi_j$}.
If $\vec x = (x_0,...,x_{j-1})^T$ is a non-zero complex vector satisfying
\begin{equation}
\label{f3_52}
(J_{5;j-1} - \lambda_0 J_{3;j-1}) \vec x = 0,
\end{equation}
for an eigenvalue $\lambda_0$, then $\vec x$ is said to be \textit{an eigenvector of the pencil $\Phi_j$
corresponding to (the eigenvalue) $\lambda_0$}.
The subspace in $\mathbb{C}^{j}$, consisting of all eigenvectors of $\Phi_j$ corresponding to $\lambda_0$ and
zero vector, is said to be \textit{the eigensubspace of the pencil $\Phi_j$
corresponding to $\lambda_0$}. 
Let us calculate $D_j(\lambda)$ in terms of $p_j,u_j$, and show that for the corresponding eigenvectors a
sort of orthogonality holds~(cf.~\cite{cit_98000_Zhedanov_JAT}).

\begin{proposition}
\label{p3_1}
Let $\Theta = \left( J_3, J_5, \alpha, \beta \right)$ be a Jacobi-type pencil, with matrices $J_3, J_5$ as in~(\ref{f1_5}),(\ref{f1_10}).
Let $\{ p_n(\lambda) \}_{n=0}^\infty$ be the associated polynomials to $\Theta$, and $\{ u_n(\lambda) \}_{n=0}^\infty$
be the shifted polynomials for $\Theta$. Then the characteristic polynomial $D_j(\lambda)$ of
the pencil $\Phi_j =J_{5;j-1} - \lambda J_{3;j-1}$ has the following form:
\begin{equation}
\label{f3_55}
D_j(\lambda) = c_j \det \left(
\begin{array}{cc}
p_j(\lambda) & u_j(\lambda)\\
p_{j+1}(\lambda) & u_{j+1}(\lambda)\end{array}
\right),\qquad c_j\in\mathbb{C},\ j\in\mathbb{N};\ (\lambda\in\mathbb{C}).
\end{equation}
The eigenspaces of the pencil $\Phi_j$ are at most two-dimensional. 
Let $\vec x'$, $\vec x''$ be some eigenvectors of the pencil $\Phi_j$, corresponding to
different eigenvalues.
Then
\begin{equation}
\label{f3_57}
\left[
J_{3;j-1} \vec x', \vec x''
\right]_C = 0,
\end{equation}
where
\begin{equation}
\label{f3_59}
\left[
\vec a, \vec b
\right]_C := (\vec a, C \vec b)_{\mathbb{C}^{j}},
\end{equation}
and $C \vec y := \overline{\vec y}$, $\vec y\in\mathbb{C}^j$.
\end{proposition}
\textbf{Proof.} 
Choose an arbitrary root $\lambda_0$ of the polynomial $D_{k+1}(\lambda)$ ($k\in\mathbb{Z}_+$). Consider 
the corresponding eigenvector $\vec x = \vec x(\lambda_0) = (x_0(\lambda_0),...,x_k(\lambda_0))$:
$$ (J_{5;k} - \lambda_0 J_{3;k}) \vec x = 0. $$
Observe that the components $x_0(\lambda_0),...,x_k(\lambda_0)$ of 
$\vec x$ satisfy the recurrent relation~(\ref{f3_5}) for $n=0,1,...,k-2$, with $\lambda=\lambda_0$. 
Since $\{ p_n(\lambda_0) \}_{n=0}^k$ and $\{ u_n(\lambda_0) \}_{n=0}^k$ satisfy the same relations, then
\begin{equation}
\label{f3_61}
\vec x(\lambda_0) = 
\widetilde c_1 
\left(
\begin{array}{ccc}
p_0(\lambda_0)\\ 
\vdots\\
p_k(\lambda_0)\end{array}
\right)
+
\widetilde c_2 
\left(
\begin{array}{ccc}
u_0(\lambda_0)\\ 
\vdots\\
u_k(\lambda_0)\end{array}
\right),\qquad
\widetilde c_1, \widetilde c_2\in\mathbb{C}.
\end{equation}
Therefore the eigensubspace corresponding to $\lambda_0$ is at most two-dimensional.
On the other hand, the components $x_0(\lambda_0),...,x_k(\lambda_0)$ of 
$\vec x$ satisfy the following two relations:
\begin{equation}
\label{f3_63}
\gamma_{k-3} x_{k-3} + (\beta_{k-2} - \lambda_0 a_{k-2}) x_{k-2} + (\alpha_{k-1} - \lambda_0 b_{k-1}) x_{k-1}
+ (\beta_{k-1} - \lambda_0 a_{k-1}) x_{k} = 0, 
\end{equation}
\begin{equation}
\label{f3_65}
\gamma_{k-2} x_{k-2} + (\beta_{k-1} - \lambda_0 a_{k-1}) x_{k-1} + (\alpha_{k} - \lambda_0 b_{k}) x_{k}
= 0. 
\end{equation}
Observe that $v_n := \widetilde c_1 p_n(\lambda_0) + \widetilde c_2 u_n(\lambda_0)$ ($n\in\mathbb{Z}_+$)
satisfy the recurrent relation~(\ref{f3_5}) with $\lambda=\lambda_0$,
for $n=k-1,k$:
$$ \gamma_{k-3} v_{k-3} + (\beta_{k-2}-\lambda_0 a_{k-2}) v_{k-2} + (\alpha_{k-1}-\lambda_0 b_{k-1}) v_{k-1} +
(\beta_{k-1}-\lambda_0 a_{k-1}) v_{k} +
$$
\begin{equation}
\label{f3_67}
+ \gamma_{k-1} v_{k+1} = 0;
\end{equation}
$$ \gamma_{k-2} v_{k-2} + (\beta_{k-1}-\lambda_0 a_{k-1}) v_{k-1} + (\alpha_k-\lambda_0 b_k) v_k +
(\beta_k-\lambda_0 a_k) v_{k+1} +
$$
\begin{equation}
\label{f3_70}
+ \gamma_k v_{k+2} = 0.
\end{equation}
Comparing relations~(\ref{f3_63}),(\ref{f3_65}) with relations~(\ref{f3_67}),(\ref{f3_70}) we conclude that
\begin{equation}
\label{f3_73}
v_{k+1} = 0,\  v_{k+2} = 0.
\end{equation}
The last relation can be rewritten in the following form:
\begin{equation}
\label{f3_75}
\left(
\begin{array}{cc}
p_{k+1}(\lambda_0) & u_{k+1}(\lambda_0)\\
p_{k+2}(\lambda_0) & u_{k+2}(\lambda_0)\end{array}
\right)
\left(
\begin{array}{cc}
\widetilde c_1\\
\widetilde c_2\end{array}
\right) = 0.
\end{equation}
Therefore
\begin{equation}
\label{f3_77}
\det \left(
\begin{array}{cc}
p_{k+1}(\lambda_0) & u_{k+1}(\lambda_0)\\
p_{k+2}(\lambda_0) & u_{k+2}(\lambda_0)\end{array}
\right)
= 0.
\end{equation}
Conversely, if relation~(\ref{f3_77}) holds with some complex $\lambda_0$, then relation~(\ref{f3_75}) is valid for some complex
$\widetilde c_1, \widetilde c_2$, not both zeros. 
Then we get relation~(\ref{f3_73}), where 
$v_n = \widetilde c_1 p_n(\lambda_0) + \widetilde c_2 u_n(\lambda_0)$, $n\in\mathbb{Z}_+$.
Therefore $\vec v := (v_0,...,v_k)^T$ satisfies the following relation:
$$ (J_{5;k} - \lambda_0 J_{3;k}) \vec v = 0,\qquad \vec v\not= 0. $$
Thus, $\lambda_0$ is a root of $D_{k+1}$.
The first statement of the proposition is proved.

%\noindent
In order to prove relation~(\ref{f3_57}) we shall proceed similar to considerations in~\cite{cit_98000_Zhedanov_JAT}.
Suppose that vectors $\vec x'$, $\vec x''$ correspond to eigenvalues $\lambda,\mu$, respectively.
We may write:
$$ \lambda \left[
J_{3;j-1} \vec x', \vec x''
\right]_C = 
\left[
J_{5;j-1} \vec x', \vec x''
\right]_C =
\left(
J_{5;j-1} \vec x', C \vec x''
\right)_{\mathbb{C}^j} = $$
$$ = 
\left(
\vec x', C C J_{5;j-1} C \vec x''
\right)_{\mathbb{C}^j} =
\left(
\vec x', C J_{5;j-1} \vec x''
\right)_{\mathbb{C}^j} =
\mu \left(
\vec x', C J_{3;j-1} \vec x''
\right)_{\mathbb{C}^j} = $$
$$ = \mu \left(
\vec x', C J_{3;j-1} C C \vec x''
\right)_{\mathbb{C}^j} =
\mu \left(
\vec x', J_{3;j-1} C \vec x''
\right)_{\mathbb{C}^j} =
\mu \left(
J_{3;j-1} \vec x', C \vec x''
\right)_{\mathbb{C}^j} = $$
$$ = \mu \left[
J_{3;j-1} \vec x', \vec x''
\right]_C. $$
Since $\lambda\not= \mu$, we conclude that $\left[
J_{3;j-1} \vec x', \vec x''
\right]_C = 0$.
$\Box$

Suppose that $\vec x$, $\vec y$ are eigenvectors of the matrix pencil $\Phi_j$, corresponding to eigenvalues $\lambda',\lambda''$, 
respectively. Then
$$ \lambda' \left(
J_{3;j-1} \vec x, \vec y
\right)_{\mathbb{C}^j} = 
\left(
J_{5;j-1} \vec x, \vec y
\right)_{\mathbb{C}^j} =
\left(
\vec x, J_{5;j-1} \vec y
\right)_{\mathbb{C}^j} = $$
\begin{equation}
\label{f3_79}
=  
\overline{\lambda''} \left(
\vec x, J_{3;j-1} \vec y
\right)_{\mathbb{C}^j} = 
\overline{\lambda''} \left( J_{3;j-1}
\vec x, \vec y
\right)_{\mathbb{C}^j}. 
\end{equation}
If $\lambda' \not= \overline{\lambda''}$ then $\left( J_{3;j-1}
\vec x, \vec y
\right)_{\mathbb{C}^j} = 0$.
In the case $J_{3;j-1} > 0$ this fact implies (by choosing $\vec x = \vec y$), as it was noticed above, that all eigenvalues are real.

\section{Non-standard three term recurrent relations and classical-type orthogonal polynomials.}

Let $\sigma$ be a non-negative measure on $\mathfrak{B}(\mathbb{R})$ with all finite power moments, 
$\int_\mathbb{R} d\sigma = 1$, $\int_\mathbb{R} |g(x)|^2 d\sigma > 0$, for any non-zero complex polynomial $g$. 
By Corollary~\ref{cc2_1}, a linear operator $\mathcal{A}$ in $L^2_\sigma$ is a model representation in $L^2_\sigma$ of the 
associated operator of some Jacobi-type pencil if and only if properties $(i)$-$(iii)$ of Theorem~\ref{tt2_1} hold.
It is readily checked that the following operator:
\begin{equation}
\label{f4_5}
\mathcal{A} [p(\lambda)] = \Lambda_0 [p(\lambda)] + p(0)[c\lambda + d],\qquad p\in\mathbb{P},
\end{equation}
where $c>-1$ and $d\in\mathbb{R}$, satisfies properties $(i)$-$(iii)$ of Theorem~\ref{tt2_1}.
By~Theorem~\ref{tt2_1} we conclude that $S(u,v)$, given by~(\ref{ff2_5}), is the spectral function of
some Jacobi-type pencil.
Such a Jacobi-type pencil $\Theta = (J_3,J_5,\alpha,\beta)$ was explicitly constructed in the proof of the Sufficiency of
Theorem~3.1 in~\cite{cit_97000_Z}.
In our case, the matrix $J_3$ is the Jacobi matrix, which corresponds to the measure $\sigma$,
and $J_5 = J_3^2$ (see the formula after~(3.8) in~\cite{cit_97000_Z}).
The constants $\alpha,\beta$ are given by the following formula:
\begin{equation}
\label{ff2_20_5}
\alpha = \frac{1}{ \xi_{0,1}\sqrt{\Delta_1} },\quad 
\beta = -\frac{\xi_{0,1} s_1 + \xi_{0,0}}{ \xi_{0,1}\sqrt{\Delta_1} }.
\end{equation}
Here $s_j$ are the power moments of $\sigma$, 
while $\Delta_n := \det (s_{k+l})_{k,l=0}^n$, $n\in\mathbb{Z}_+$; $\Delta_{-1}:=1$, are the Hankel determinants.
The coefficients $\xi_{k,j}$ are taken from property~(ii) of Theorem~\ref{tt2_1}.

\begin{theorem}
\label{t4_1}
Let $\sigma$ be a non-negative measure on $\mathfrak{B}(\mathbb{R})$ with all finite power moments, 
$\int_\mathbb{R} d\sigma = 1$, $\int_\mathbb{R} |g(x)|^2 d\sigma > 0$, for any non-zero complex polynomial $g$;
and $\mathcal{A}$ be given by~(\ref{f4_5}). 
Define $S(u,v)$ by relation~(\ref{ff2_5}). 
Let $\Theta = (J_3,J_5,\alpha,\beta)$ be the Jacobi-type pencil with the spectral function $S$, constructed by~(\ref{f4_5})-(\ref{ff2_20_5}).
Denote by $\{ p_n(\lambda) \}_{n=0}^\infty$ the associated polynomials to the pencil $\Theta$,
and denote by $\{ r_n(\lambda) \}_{n=0}^\infty$ the orthonormal polynomials (with positive leading coefficients) with
respect to the measure $\sigma$. Then
\begin{equation}
\label{f4_7}
p_n(\lambda) = \frac{1}{c+1} r_n(\lambda) - \frac{d}{c+1} \frac{r_n(\lambda) - r_n(0)}{\lambda} +
\frac{c}{c+1} r_n(0),\qquad n\in\mathbb{Z}_+;
\end{equation}
\begin{equation}
\label{f4_9}
r_n(\lambda) = (c+1) p_n(\lambda) + (c+1)d \frac{p_n(\lambda) - p_n(d)}{\lambda - d} -
c p_n(d),\qquad n\in\mathbb{Z}_+.
\end{equation}
In~(\ref{f4_7}),(\ref{f4_9}) we mean the limit expressions at $\lambda = 0$ and $\lambda = d$, respectively.
The following recurrent relation, involving three subsequent associated polynomials, holds:
\begin{equation}
\label{f4_15}
\lambda p_n(\lambda) = \frac{ p_n(d) }{c+1} (c\lambda + d) + a_{n-1} p_{n-1}(\lambda) + b_n p_n(\lambda) + a_n p_{n+1}(\lambda),\ 
n\in\mathbb{Z}_+,\ (\lambda\in\mathbb{C}).
\end{equation}
The following orthonormality relations hold:
$$ \int_{\mathbb{R}\backslash\{ d \}} (p_n(\lambda),p_n(d))
\left(
\begin{array}{cc} (c+1)^2 \left( \frac{\lambda}{\lambda - d} \right)^2 & (-c-1) \frac{\lambda (c\lambda + d)}{(\lambda - d)^2} \\
(-c-1) \frac{\lambda (c\lambda + d)}{(\lambda - d)^2} & \left( \frac{ c\lambda + d }{\lambda - d} \right)^2 
\end{array}
\right)
\left(
\begin{array}{cc} p_m(\lambda) \\
p_m(d)
\end{array}
\right) d\sigma + $$
\begin{equation}
\label{f4_18}
+ 
(p_n(d), p_n'(d))
\left(
\begin{array}{cc} 1 & (c+1) d \\
(c+1) d  & (c+1)^2 d^2 
\end{array}
\right)
\left(
\begin{array}{cc} p_m(d) \\
p_m'(d)
\end{array}
\right) \sigma(\{ d \}) = \delta_{n,m},\ n,m\in\mathbb{Z}_+.
\end{equation}
\end{theorem}
\textbf{Proof.} 
In view of relation~(\ref{f3_160}) it will be useful to calculate $p_n(\mathcal{A})(1)$.
By the induction argument one can verify that
\begin{equation}
\label{f4_25}
\mathcal{A}^n [1] =
\left[
(c+1) 
\left\{
\lambda^n + d\lambda^{n-1} + d^2\lambda^{n-2} + ... + d^{n-1}\lambda
\right\}
+ d^n
\right],\qquad n\in\mathbb{N}.
\end{equation}
Relation~(\ref{f4_25}) can be written in a more compact form (and including the case $n=0$):
\begin{equation}
\label{f4_27}
\mathcal{A}^n [1] =
\left[
(c+1) 
\frac{ \lambda (d^n - \lambda^n) }{d-\lambda}
+ d^n
\right],\qquad n\in\mathbb{Z}_+.
\end{equation}
Here we mean the limit expression for $\lambda = d$.
Let $u(\lambda) = \sum_{k=0}^r a_k \lambda^k$ be an arbitrary complex polynomial of degree $r$.
Then
$$ u(\mathcal{A})[1] = \sum_{k=0}^r a_k \mathcal{A}^k [1] = 
\left[
\frac{(c+1)\lambda}{ d-\lambda } \sum_{k=0}^r (a_k d^k - a_k \lambda^k)
+ \sum_{k=0}^r a_k d^k
\right] =
$$
$$ =
\left[
\frac{(c+1)\lambda}{ d-\lambda } (u(d) - u(\lambda)) + u(d)
\right] = $$
\begin{equation}
\label{f4_29}
= \left[ (c+1) u(\lambda) + (c+1) d \frac{ u(\lambda) - u(d) }{\lambda - d} - c u(d)
\right]. 
\end{equation}
By relation~(\ref{f3_160}) we conclude that the following orthogonality relations hold:
$$ \int_{\mathbb{R}} 
\left\{
(c+1) p_n(\lambda) + (c+1) d \frac{ p_n(\lambda) - p_n(d) }{\lambda - d} - c p_n(d)
\right\} * $$
\begin{equation}
\label{f4_30}
*
\left\{
(c+1) p_m(\lambda) + (c+1) d \frac{ p_m(\lambda) - p_m(d) }{\lambda - d} - c p_m(d)
\right\}
d\sigma
= \delta_{n,m},\ n,m\in\mathbb{Z}_+.
\end{equation}
Observe that expressions in brackets $\{ ... \}$ in~(\ref{f4_30}) are real polynomials of degrees $n$ and $m$ with
positive leading coefficients.
Since the orthonormal system of polynomials with respect to $\sigma$ (having positive leading coefficients) is unique
(e.g.~\cite{cit_20000_Suetin}) then relation~(\ref{f4_9}) holds.

Multiplying relation~(\ref{f4_9}) by $\lambda - d$ and simplifying the resulting expression we get
\begin{equation}
\label{f4_32}
(c+1) \lambda p_n(\lambda) - (c\lambda + d) p_n(d) = r_n(\lambda) (\lambda - d),\qquad n\in\mathbb{Z}_+.
\end{equation}
For $\lambda = 0$ we obtain that $-d p_n(d) = -d r_n(0)$. If $d\not=0$, then $p_n(d) = r_n(0)$.
If $d=0$, then relation~(\ref{f4_32}) becomes 
$$ (c+1) \lambda p_n(\lambda) - c\lambda p_n(0) = \lambda r_n(\lambda),\qquad n\in\mathbb{Z}_+. $$
For $\lambda\not= 0$ we may write: 
\begin{equation}
\label{f4_34}
(c+1) p_n(\lambda) - c p_n(0) = r_n(\lambda),\qquad n\in\mathbb{Z}_+.
\end{equation}
Passing to the limit as $\lambda\rightarrow 0$, we obtain that $p_n(0) = r_n(0)$.
Thus, for arbitrary $d\in\mathbb{R}$ it holds that
\begin{equation}
\label{f4_36}
p_n(d) = r_n(0),\qquad n\in\mathbb{Z}_+.
\end{equation}
By relations~(\ref{f4_36}),(\ref{f4_32}) we obtain that relation~(\ref{f4_7}) is valid.

Multiplying the three-term recurrent relation for $r_n(\lambda)$ by $(\lambda-d)$ and using relation~(\ref{f4_32}),
we obtain that
$$ (c+1) \lambda (a_{n-1} p_{n-1}(\lambda) + (b_n - \lambda) p_n(\lambda) + a_n p_{n+1}(\lambda)) = $$
$$ = (c\lambda + d) (a_{n-1} p_{n-1}(d) + (b_n - \lambda) p_n(d) + a_n p_{n+1}(d)),\qquad n\in\mathbb{Z}_+. $$
By relation~(\ref{f4_36}) and the recurrent relation for $r_n$ we get:
$$ (c+1) \lambda (a_{n-1} p_{n-1}(\lambda) + (b_n - \lambda) p_n(\lambda) + a_n p_{n+1}(\lambda)) = $$
$$ = (c\lambda + d) (-\lambda) p_n(d),\qquad n\in\mathbb{Z}_+. $$
Therefore recurrent relation~(\ref{f4_15}) is valid.            
Orthogonality relations~(\ref{f4_18}) are equivalent to relations~(\ref{f4_30}), since expressions under the integrals
coincide (the term for $\lambda = d$ is written separately).
$\Box$

The following example shows that polynomials $p_n(\lambda)$ from Theorem~\ref{t4_1} can have multiple or complex
roots.

\begin{example}
\label{e4_1}
Consider the Jacobi-type pencil $\Theta$ from Theorem~\ref{t4_1} with an additional assumption $c=0$.
The recurrent relation~(\ref{f4_15}) takes the following form:
\begin{equation}
\label{f4_38}
\lambda p_n(\lambda) = d p_n(d) + a_{n-1} p_{n-1}(\lambda) + b_n p_n(\lambda) + a_n p_{n+1}(\lambda),\ 
n\in\mathbb{Z}_+,\ (\lambda\in\mathbb{C}).
\end{equation}
Using this relation we calculate several first polynomials $p_n(\lambda)$:
$$ p_0(\lambda) = 1,\quad p_1(\lambda) = \frac{1}{a_0} (\lambda - b_0 - d), $$
\begin{equation}
\label{f4_42}
p_2(\lambda) = \frac{1}{a_0 a_1} \left\{
\lambda^2 - (b_0 + b_1 + d) \lambda + (b_0 + b_1) d + b_0 b_1 - a_0^2
\right\}.
\end{equation}
The discriminant $D$ of the quadratic expression in the brackets in~(\ref{f4_42}) is equal to
\begin{equation}
\label{f4_46}
D = d (d - 2b_0 - 2b_1) + (b_0 - b_1)^2 + 4 a_0^2.
\end{equation}
If the following assumptions hold:
$$ a_0 = \frac{1}{2},\ b_0 = b_1 = 1,\ d=1, $$
then $D=-2$.
On the other hand, if
$$ a_0 = \frac{\sqrt{3}}{2},\ b_0 = b_1 = 1,\ d=1, $$
then $D=0$.
Thus, the polynomial $p_2$ can have multiple or complex roots.
\end{example}

Consider the differential equation~(\ref{f1_150}), where $\lambda$ is a real parameter and
$$ c(t) = \sum_{j=0}^4 c_j t^j,\quad d(t) = \sum_{j=0}^3 d_j t^j,\quad f(t) = \sum_{j=0}^2 f_j t^j, $$
\begin{equation}
\label{f4_48}
g(t) = g_1 t + g_0,\quad h(t) = h_0;
\end{equation}
\begin{equation}
\label{f4_49}
\varphi(t) = \varphi_2 t^2 + \varphi_1 t + \varphi_0,\quad \psi(t) = \psi_1 t + \psi_0,\quad \theta(t) = \theta_0,
\end{equation}
where all $c_j,d_j,f_j,g_j,h_j,\varphi_j,\psi_j,\theta_j$ are some real numbers.
The following lemma, which reflects the idea from~\cite{cit_5_Azad}, will be useful.

\begin{lemma}
\label{l4_1}
Let the differential equation~(\ref{f1_150}) be given, with a real parameter $\lambda$ and
some polynomial coefficients of the form~(\ref{f4_48}),(\ref{f4_49}).
Let $n$ be an arbitrary non-negative integer.
A polynomial of the following form:
\begin{equation}
\label{f4_51}
y(t) = \sum_{k=0}^n \mu_k t^k,\qquad \mu_k\in\mathbb{R},
\end{equation}
is a solution of equation~(\ref{f1_150})
if and only if
$\{ \mu_k \}_{k=0}^n$ is a solution of the following system of linear equations:
$$ 
\{
j(j-1)(j-2)(j-3) c_4 + j(j-1)(j-2) d_3 + j(j-1) f_2 + j g_1 + h_0 
+ $$
$$ + \lambda
(
j(j-1) \varphi_2 + j \psi_1 + \theta_0
)
\}
\mu_j +
$$
$$ 
\{
(j+1)j(j-1)(j-2) c_3 + (j+1)j(j-1) d_2 + (j+1)j f_1 + (j+1) g_0 
+ $$
$$ + \lambda
(
(j+1)j \varphi_1 + (j+1) \psi_0
)
\}
\mu_{j+1} +
$$
$$ 
\{
(j+2)(j+1)j(j-1) c_2 + (j+2)(j+1)j d_1 + (j+2)(j+1) f_0
+ $$
$$ + \lambda
(j+2)(j+1) \varphi_0
\}
\mu_{j+2} +
$$
$$ 
\{
(j+3)(j+2)(j+1)j c_1 + (j+3)(j+2)(j+1) d_0
\}
\mu_{j+3} +
$$
\begin{equation}
\label{f4_55}
+ (j+4)(j+3)(j+2)(j+1) c_0 \mu_{j+4} = 0,\qquad j=0,1,...,n.
\end{equation}
Here $\mu_{n+1} = \mu_{n+2} = \mu_{n+3} = 0$. 

\end{lemma}
\textbf{Proof.}
It is enough to substitute the expression for $y(t)$ into the left-hand side of differential equation~(\ref{f1_150})
and calculate the coefficients by $t^j$, $j=0,1,...,n$. The coefficient by $t^j$ equals to the left-hand side
of~(\ref{f4_55}).
$\Box$

Consider the Jacobi-type pencil $\Theta$ from Theorem~\ref{t4_1} with an additional assumption $d=0$.
Relations~(\ref{f4_7}),(\ref{f4_9}) take the following form:
\begin{equation}
\label{f4_56}
p_n(\lambda) = \frac{1}{c+1} r_n(\lambda) +
\frac{c}{c+1} r_n(0),\qquad n\in\mathbb{Z}_+;
\end{equation}
\begin{equation}
\label{f4_57}
r_n(\lambda) = (c+1) p_n(\lambda) - c p_n(0),\qquad n\in\mathbb{Z}_+.
\end{equation}
Suppose that $r_n$ satisfy the following differential equation:
\begin{equation}
\label{f4_59}
v (t) r_n''(t) + w (t) r_n'(t) + \gamma_n r_n(t) = 0,\qquad n\in\mathbb{Z}_+,
\end{equation}
where $v(t), w(t)$ are some real polynomials, $\gamma_n\in\mathbb{R}$. 
After the differentiation we obtain that $r_n$ satisfy the third order differential equation:
\begin{equation}
\label{f4_61}
v (t) y'''(t) + (v'(t) + w (t)) y''(t) + w'(t) y'(t) + \gamma_n y'(t) = 0,\qquad n\in\mathbb{Z}_+.
\end{equation}
Since the derivatives of $r_n$ and $p_n$ of orders $1,2,...,$ differ by a constant factor, then $p_n$
satisfy equation~(\ref{f4_61}).
Thus, all classical orthogonal polynomials $r_n$ admit generalizations $p_n$ by formula~(\ref{f4_56}),
which satisfy the third order differential equation~(\ref{f4_61}).
If $c\not= 0$, then $p_n$ is not a unique (up to a constant multiple) real polynomial solution of order $n$ for differential equation~(\ref{f4_61}).
There exists a case where the unicity (up to a constant multiple) takes place. It is described in the following theorem.

\begin{theorem}
\label{t4_2}
In assumptions of Theorem~\ref{t4_1} we additionally suppose that $\sigma$ and $J_3$ correspond to
orthonormal Jacobi polynomials $r_n(\lambda) = P_n(\lambda;a,b)$ ($a,b>-1$)
and $c=0$; $d=1$.
In this case, the associated polynomial $p_n$ ($n\in\mathbb{Z}_+$):
\begin{equation}
\label{f4_65}
p_n(\lambda) = r_n(\lambda) - \frac{r_n(\lambda) - r_n(0)}{\lambda},
\end{equation}
is a unique, up to a constant multiple, real $n$-th order polynomial solution of the following $4$-th order differential equation:
$$ -(t+1)t(t-1)^2 y^{(4)}(t) + (t-1)(-(a+b+10)t^2 + (b-a)t + 4) y^{(3)}(t) + $$
$$ + (-3(2a+2b+8)t^2 + (a+9b+22)t + 3a-3b) y''(t) + (-6(a+b+2)t + 2a+6b + 8) y'(t) + $$
\begin{equation}
\label{f4_67}
+ \lambda_n 
(
t(t-1) y''(t) + 2 (2t-1) y'(t) + 2 y(t)
) = 0,
\end{equation}
where $\lambda_n = n(n+a+b+1)$.

Moreover, there exists a unique $\lambda_n\in\mathbb{R}$, such that differential equation~(\ref{f4_67})
has a real $n$-th order polynomial solution. 
\end{theorem}
\textbf{Proof.} 
The orthonormal Jacobi polynomials $r_n(t)$ satisfy the following differential equation:
\begin{equation}
\label{f4_70}
(1-t^2) y''(t) + (b-a-(a+b+2)t) y'(t) + n(n+a+b+1) y(t) = 0.
\end{equation}
Observe that
$$ r_n(t) = \frac{tp_n(t) - p_n(1)}{t-1},\ 
r_n'(t) = (p_n(t) + t p_n'(t)) \frac{1}{t-1} - (tp_n(t) - p_n(1)) \frac{1}{(t-1)^2},
$$
$$ (1-t^2) r_n''(t) = 
-(t+1)
\left\{
2 p_n'(t) + t p_n''(t) - \frac{2( p_n(t) + t p_n'(t) )}{t-1} + \frac{ 2(tp_n(t) - p_n(1)) }{(t-1)^2}
\right\}. $$
Substituting the latter expressions into equation~(\ref{f4_70}), multiplying by $(t-1)^2$ and simplifying we get
$$ -(t+1) t (t-1)^2 p_n''(t) + $$
$$+ (
-2(t+1)(t-1)^2 + 2(t+1)(t-1)t + (b-a-(a+b+2)t) t(t-1)
) p_n'(t) + $$
$$ + (
2(t+1)(t-1) - 2 (t+1)t - (b-a-(a+b+2)t)
) p_n(t)
+ n(n+a+b+1) t(t-1) p_n(t) + $$
$$ + d(t), $$
where $d(t) := ( 2(t+1) + (b-a-(a+b+2)t - n(n+a+b+1)(t-1)) p_n(1)$.
Differentiating the latter expression two times we can remove $d(t)$ and obtain equation~(\ref{f4_67}).

In order to apply Lemma~\ref{l4_1} we notice that in our case it holds:  
$$ c_4 = -1,\ c_3 = 1,\ c_2 = 1,\ c_1 = -1,\ c_0 = 0; $$
$$ d_3 = - (a+b+10),\ d_2 = 2b+10,\ d_1 = a-b+4,\ d_0 = -4; $$
$$ f_2 = -6(a+b+4),\ f_1 = a+9b+22,\ f_0 = 3a -3b; $$
$$ g_1 = -6(a+b+2),\ g_0 = 2a+6b+8;\ h_0 = 0; $$
$$ \varphi_2 = 1,\ \varphi_1 = -1,\ \varphi_0 = 0;\ \psi_1 = 4,\ \psi_0 = -2;\ \theta_0 = 2. $$
Equation~(\ref{f4_55}) for $j=n$ uniquely determines $\lambda = n(n+a+b+1)$.
Equations~(\ref{f4_55}) with $j=n-1,n-2,...,0$, uniquely determine the coefficients of a polynomial solution of equation~(\ref{f4_67}).
Thus, $p_n(t)$ is a unique real $n$-th order polynomial solution of equation~(\ref{f4_67}).
$\Box$

\section{High-order difference equations and matrix orthogonal polynomials.}

Let $\{ p_n(\lambda) \}_{n=0}^\infty$ ($p_n$ has degree $n$ and a positive leading coefficient)
is a set of complex polynomials satisfying the following difference equation:
\begin{equation}
\label{f5_7}
\sum_{j=1}^N (\overline{\alpha_{k-j,j}} p_{k-j} (\lambda) + \alpha_{k,j} p_{k+j}(\lambda))
+ \alpha_{k,0} p_k(\lambda) = \lambda^N p_k(\lambda), 
\end{equation}
with some complex coefficients $\alpha_{m,n}$ ($m\in\mathbb{Z}_+$, $n=0,1,...,N$): $\alpha_{m,N} > 0$, $\alpha_{m,0}\in\mathbb{R}$;
$N\in\mathbb{N}$. Here $\alpha_{m,n}$ and $p_k$ with negative indices are zero. 
Equation~(\ref{f5_7}) can be written in the following matrix form:
\begin{equation}
\label{f5_9}
J \vec p(\lambda) = \lambda^N \vec p(\lambda),\qquad \vec p(\lambda) =
\left(
\begin{array}{ccc}
p_0(\lambda)\\
p_1(\lambda)\\
\vdots
\end{array}
\right),
\end{equation}
where $J$ is a complex Hermitian $(2N+1)$-diagonal semi-infinite matrix.
Observe that $J$ can be viewed as a block Jacobi matrix. 
Thus, $\{ p_n(\lambda) \}_{n=0}^\infty$ have a transparent algebraic connection with the corresponding
orthonormal matrix polynomials $\{ P_n(x) \}_{n=0}^\infty$~\cite{cit_98300_D_A}:
\begin{equation}
\label{f5_12}
P_n(x) =
\left(
\begin{array}{ccc}
R_{N,0}(p_{nN})(x) & \ldots & R_{N,N-1}(p_{nN})(x)\\
\vdots & \ldots & \vdots \\
R_{N,0}(p_{nN+N-1})(x) & \ldots & R_{N,N-1}(p_{nN+N-1})(x)
\end{array}
\right),
\end{equation}
where 
$$ R_{N,m}(p)(t) := \sum_n \frac{ p^{(nN+m)}(0) }{ (nN+m)! } t^n,\qquad p\in\mathbb{P}. $$
Denote $P_N = \{ \lambda\in\mathbb{C}:\ \lambda^N\in\mathbb{R} \}$, and
let $\varepsilon$ be a primitive $N$-th root of unity.
Polynomials $\{ p_n(\lambda) \}_{n=0}^\infty$ are orthonormal with respect to the following functional:
$$ B(u,v) = \int_{P_N}
(u(\lambda), u(\lambda\varepsilon), \ldots, u(\lambda\varepsilon^{N-1})) dW(\lambda)
\overline{
\left(
\begin{array}{cccc}
v(\lambda) \\
v(\lambda\varepsilon) \\
\vdots \\
v(\lambda\varepsilon^{N-1})
\end{array}
\right)} + $$
\begin{equation}
\label{f5_14}
+ 
(u(0), u'(0), \ldots, u^{(N-1)}(0)) M
\overline{
\left(
\begin{array}{cccc}
v(0) \\
v'(0) \\
\vdots \\
v^{(N-1)}(0)
\end{array}
\right)},\qquad
p,q\in\mathbb{P},
\end{equation}
where $W(\lambda)$ is a non-decreasing matrix-valued function on $P_N\backslash\{ 0 \}$,
$M\geq 0$ is a $(N\times N)$ complex matrix (at $\lambda = 0$ the integral is understood as improper).
The sesquilinear functional $B$ has the following properties: $\overline{B(u,v)} = B(v,u)$, and
\begin{equation}
\label{f5_16}
B(\lambda^N u(\lambda), v(\lambda)) = B(u(\lambda), \lambda^N v(\lambda)),\qquad u,v\in\mathbb{P}.
\end{equation}

\begin{theorem}
\label{t5_1}
Let $\Theta = (J_3,J_5,\alpha,\beta)$ be a Jacobi-type pencil. 
Denote by $\{ p_n(\lambda) \}_{n=0}^\infty$ the associated polynomials to the pencil $\Theta$,
and denote by $\mathcal{A}=\mathcal{A}_\sigma$ the associated operator for $\Theta$.
Polynomials $\{ p_n(\lambda) \}_{n=0}^\infty$ satisfy recurrent relation~(\ref{f5_7}) with a positive integer $N$
and with some complex coefficients $\alpha_{m,n}$ ($\alpha_{m,N} > 0$, $\alpha_{m,0}\in\mathbb{R}$)
if and only if $\mathcal{A}^N$ is a symmetric operator.
\end{theorem}
\textbf{Proof.} 
\textit{Necessity.}
Suppose that the associated polynomials $p_n(\lambda)$ satisfy recurrent relation~(\ref{f5_7}), with a positive integer $N$
and with some complex coefficients $\alpha_{m,n}$.
Denote by $S(u,v)$ the spectral function of $\Theta$.
Since $p_n$ are orthonormal with respect to sesquilinear functionals $B$ (from~(\ref{f5_14})) and $S$, then
$B=S$.
By~(\ref{f5_16}),(\ref{ff2_5}) we may write:
\begin{equation}
\label{f5_18}
\int_{\mathbb{R}} (\lambda^N u(\lambda))(\mathcal{A})(1) \overline{ v(\mathcal{A})(1) } d\sigma
=
\int_{\mathbb{R}} u(\mathcal{A})(1) \overline{ (\lambda^N v(\lambda))(\mathcal{A})(1) } d\sigma, 
\end{equation}
for arbitrary $u,v\in\mathbb{P}$.
Since
\begin{equation}
\label{f5_19}
(\lambda^N u(\lambda))(\mathcal{A})[1] =
\mathcal{A}^N \left\{ u(\mathcal{A})[1] \right\},\qquad u\in\mathbb{P},
\end{equation}
then
\begin{equation}
\label{f5_20}
\int_{\mathbb{R}} \mathcal{A}^N \left\{
u(\mathcal{A})(1)
\right\}
\overline{ v(\mathcal{A})(1) } d\sigma
=
\int_{\mathbb{R}} u(\mathcal{A})(1) \overline{ 
\mathcal{A}^N 
\left\{
v(\mathcal{A})(1)
\right\}
} d\sigma,\quad u,v\in\mathbb{P}.  
\end{equation}
It remains to check that for an arbitrary complex polynomial $w(\lambda)$ there exists a representation: $[w] = u(\mathcal{A})[1]$ with a
suitable $u\in\mathbb{P}$.
Observe that polynomials
$\widetilde g_n$:
$$ [\widetilde g_n] = \mathcal{A}^n [1],\qquad n\in\mathbb{Z}_+, $$
form a linear basis in $\mathbb{P}$, see the formula following~(3.7) in~\cite{cit_97000_Z}.
Then
$$ [w] = \sum_{j} \widetilde a_j [\widetilde g_j] = \left(
\sum_j \widetilde a_j \mathcal{A}^j
\right) [1], $$
for some complex numbers $\widetilde a_j$.

\noindent
\textit{Sufficiency.}
Suppose that $\mathcal{A}^N$ is symmetric ($n\in\mathbb{N}$). Then
\begin{equation}
\label{f5_22}
\int_\mathbb{R}
\mathcal{A}^N [w(\lambda)]
\overline{[g(\lambda)]}
d\sigma
=
\int_\mathbb{R}
[w(\lambda)]
\overline{ \mathcal{A}^N [g(\lambda)] }
d\sigma,\quad w,g\in\mathbb{P}. 
\end{equation}
In particular, relation~(\ref{f5_22}) implies relation~(\ref{f5_20}).
By~(\ref{f5_19}) we conclude that relation~(\ref{f5_18}) holds.
Relation~(\ref{f5_18}) means that the spectral function $S$ satisfies the following relation:
\begin{equation}
\label{f5_25}
S(\lambda^N u(\lambda), v(\lambda)) = S(u(\lambda), \lambda^N v(\lambda)),\qquad u,v\in\mathbb{P}.
\end{equation}
%Using this relation, it is not hard to obtain recurrent relation~(\ref{f5_7}). 
Since $p_k$ are real polynomials, $\deg p_k = k$, we may write:
\begin{equation}
\label{f5_25_3}
\lambda^N p_k(\lambda) = \sum_{i=0}^{k+N} \xi_{k,i} p_i(\lambda),\qquad \xi_{k,i}\in\mathbb{R}.
\end{equation}
By the orthogonality we get
\begin{equation}
\label{f5_25_5}
\xi_{k,j} = S(\lambda^N p_k(\lambda), p_j(\lambda)),\qquad 0\leq j\leq k+N.
\end{equation}
For $i$: $0\leq i< k-N$, we may write:
$$ \xi_{k,i} = S(\lambda^N p_k(\lambda), p_i(\lambda)) = S(p_k(\lambda), \lambda^N p_i(\lambda)) = 0. $$
Then
\begin{equation}
\label{f5_25_7}
\lambda^N p_k(\lambda) = \sum_{i=k-N}^{k+N} \xi_{k,i} p_i(\lambda)
= \sum_{j=1}^N (\xi_{k,k-j} p_{k-j}(\lambda) + \xi_{k,k+j} p_{k+j}(\lambda)) + \xi_{k,k} p_{k}(\lambda).
\end{equation}
Set 
\begin{equation}
\label{f5_25_9}
\alpha_{k,0} = \xi_{k,k},\quad \alpha_{k,j} = \xi_{k,k+j},\qquad j=1,2,...,N;\ (k\in\mathbb{Z}_+). 
\end{equation}
Observe that
$$ \xi_{k,k-j} = S(\lambda^N p_k(\lambda), p_{k-j}(\lambda)) = S(p_k(\lambda), \lambda^N p_{k-j}(\lambda)) = $$
\begin{equation}
\label{f5_25_11}
= \overline{ S(\lambda^N p_{k-j}(\lambda), p_k(\lambda)) } = \overline{ \xi_{k-j,k} } = \overline{\alpha_{k-j,j}}, 
\end{equation}
for $j=1,...,N$: $k-j \geq 0$.  
By~(\ref{f5_25_7})-(\ref{f5_25_11}) we conclude that relation~(\ref{f5_7}) holds.
$\Box$

Let us return to polynomials $p_n$ from the previous section. Namely, assume that
assumptions of Theorem~\ref{t4_1} hold. Let $\mathcal{A}=\mathcal{A}_\sigma$ be the associated operator for the pencil $\Theta$.
Assume additionally that $c\not=0$ and $d\not=0$.
Let us show that for an arbitrary $N\in\mathbb{N}$ the operator $\mathcal{A}^N$ is not symmetric. 
In fact, by the induction argument one can check that
\begin{equation}
\label{f5_27}
\mathcal{A}^n [p(\lambda)] = 
\left[
\lambda^n p(\lambda) + p(0) (c\lambda + d)
( \lambda^{n-1} + \lambda^{n-2} d + ... + d^{n-1} )
\right],\ n\in\mathbb{N}.
\end{equation}
Relation~(\ref{f5_27}) can be written in a compact form:
\begin{equation}
\label{f5_29}
\mathcal{A}^n [p(\lambda)] = 
\Lambda_0^n [p(\lambda)] +
p(0) 
\left[
(c\lambda + d)
\frac{ d^n - \lambda^n }{ d-\lambda }
\right],\qquad n\in\mathbb{N}.
\end{equation}
Thus, $\mathcal{A}^N$ is not symmetric if the following operator
$$ B [p] := 
p(0) 
\left[
(c\lambda + d)
\frac{ d^N - \lambda^N }{ d-\lambda }
\right],\qquad (p\in\mathbb{P}),  $$
is not symmetric.
Denote
$$ \varphi(\lambda) = (c\lambda + d)
\frac{ d^N - \lambda^N }{ d-\lambda }. $$
Let $\{ r_n(x) \}_0^\infty$ ($r_0=1$) be orthonormal polynomials (having positive leading coefficients)
with respect to $\sigma$ (where $\sigma$, as usual, corresponds to $J_3$).
Choose an arbitrary $m>N$ such that $r_m(0)\not= 0$. It is always possible according to the three-term recurrent relation for $r_n$.
Then
$$ (B [\varphi], [r_m])_{L^2_\sigma} = d^N ([\varphi], [r_m])_{L^2_\sigma} = 0, $$
since $\varphi$ is a polynomial of degree at most $N$.
On the other hand,
$$ ([\varphi], B [r_m])_{L^2_\sigma} = r_m(0) ([\varphi], [\varphi])_{L^2_\sigma} \not= 0. $$
Thus, $\mathcal{A}^N$ is not symmetric.

\noindent
\textbf{Acknowledgements.} The author is grateful to Prof. Zhedanov for pointing him references~\cite{cit_1000_Chaundy} 
and~\cite{cit_15000_Sawyer}.

\begin{center}
{\large\bf 
Difference equations related to Jacobi-type pencils.}
\end{center}
\begin{center}
{\bf S.M. Zagorodnyuk}
\end{center}

In this paper we study various difference equations related to Jacobi-type pencils. By a Jacobi-type pencil one means
the following pencil:
$J_5 - \lambda J_3$,
where 
$J_3$ is a Jacobi matrix and
$J_5$ is a semi-infinite real symmetric five-diagonal matrix with positive numbers on the second subdiagonal.
The basic set of solutions for the corresponding $4$-th order difference equation is constructed. 
Spectral properties of the truncated pencil and some special matrix orthogonality relations are investigated.
Classical type orthogonal polynomials satisfying a $4$-th order differential equation are constructed.

}

%%%%%%%%%%%%%%%%%%%%%%%%%%%     23:58   7.02.2018


\begin{thebibliography}{99}

%\bibitem{cit_5_Atkinson}
%Atkinson, F. V. Discrete and Continuous Boundary Problems. 
%Academic Press, New York, London, 1964. (Russian edition: Mir, Moscow, 1968).

\bibitem{cit_5_Azad}
Azad, H.; Laradji, A.; Mustafa, M. T. Polynomial solutions of differential equations. Adv. Difference Equ. 2011, 2011:58, 12 pp.

\bibitem{cit_8_Ben_Amara_Vladimirov_Shkalikov}
Ben Amara J., Vladimirov A. A., Shkalikov A. A. 
Spectral and Oscillatory Properties of a Linear Pencil of Fourth-Order Differential Operators. Math. Notes, 94:1 (2013), 49--59. 

%\bibitem{cit_10_Berezansky_Book}
%Berezansky, Yu. M. Expansions in Eigenfunctions of Selfadjoint Oper-
%ators. Amer. Math. Soc., Providence, RI, 1968. (Russian edition: Naukova
%Dumka, Kiev, 1965).

\bibitem{cit_1000_Chaundy}
Chaundy, T. W. Second-order linear differential equations with polynomial solutions. Quart. J. Math., Oxford Ser. (2) 4, (1953). 81--95.


\bibitem{cit_3000_Chihara}
Chihara, T. S. An introduction to orthogonal polynomials. Mathematics and its Applications, Vol. 13. 
Gordon and Breach Science Publishers, New York-London-Paris, 1978. xii+249 pp.

\bibitem{cit_98100_Dereviagin}
Derevyagin, Maxim; Tsujimoto, Satoshi; Vinet, Luc; Zhedanov, Alexei. 
Bannai-Ito polynomials and dressing chains. Proc. Amer. Math. Soc. 142 (2014), no. 12, 4191--4206.

\bibitem{cit_98300_D_A}
Dur\'an, A. J.; Van Assche, W. Orthogonal matrix polynomials and higher-order recurrence relations. 
Linear Algebra Appl. 219 (1995), 261--280.


%\bibitem{cit_3500_Ch_Z}
%Choque Rivero, Abdon E.; Zagorodnyuk, Sergey M. Orthogonal polynomials on rays: Christoffel's formula. 
%Bol. Soc. Mat. Mexicana (3) 15 (2009), no. 2, 149--164.

%\bibitem{cit_4000_D_P_S}
%Damanik, David; Pushnitski, Alexander; Simon, Barry. The analytic theory of matrix orthogonal polynomials. 
%Surv. Approx. Theory 4 (2008), 1--85.

\bibitem{cit_5000_Ismail}
Ismail, Mourad E. H. Classical and quantum orthogonal polynomials in one variable. With two chapters by Walter Van Assche. 
With a foreword by Richard A. Askey. Encyclopedia of Mathematics and its Applications, 98. Cambridge University Press, Cambridge, 2005. xviii+706 pp.

\bibitem{cit_5100_Krall}
Krall, Allan M. Hilbert space, boundary value problems and orthogonal polynomials. 
Operator Theory: Advances and Applications, 133. Birkhäuser Verlag, Basel, 2002. xiv+352 pp.

%\bibitem{cit_7000_Markus}
%Markus, A. S. Introduction to the spectral theory of polynomial operator pencils. Translated from the Russian by H. H. McFaden. 
%Translation edited by Ben Silver. With an appendix by M. V. Keldysh. Translations of Mathematical Monographs, 71. American Mathematical Society, 
%Providence, RI, 1988. iv+250 pp.

%\bibitem{cit_7100_M_G_A}
%Messirdi B., Gherbi A., Amouch M. A spectral analysis of linear operator pencils on Banach spaces with application to
%quotient of bounded operators. International Journal of Analysis and Applications, {\bf 7}, 2 (2015), 104--128.



%\bibitem{cit_Piv}
%M\"oller, Manfred; Pivovarchik, Vyacheslav. 
%Spectral theory of operator pencils, Hermite-Biehler functions, and their applications. 
%Operator Theory: Advances and Applications, 246. Birkh\"auser/Springer, Cham, 2015. xvii+412 pp.

%\bibitem{cit_7700_Parlett}
%Parlett, Beresford N. The symmetric eigenvalue problem. Corrected reprint of the 1980 original. 
%Classics in Applied Mathematics, 20. Society for Industrial and Applied Mathematics (SIAM), Philadelphia, PA, 1998. {\rm xxiv}+398 pp.


%\bibitem{cit_8000_Rodman}
%Rodman, Leiba. An introduction to operator polynomials. 
%Operator Theory: Advances and Applications, 38. Birkhäuser Verlag, Basel, 1989. {\rm xii}+389 pp. 



%\bibitem{cit_10000_Simon}
%Simon, Barry. Szeg\"o's theorem and its descendants. Spectral theory for $L\sp 2$ perturbations of orthogonal polynomials. 
%M. B. Porter Lectures. Princeton University Press, Princeton, NJ, 2011. xii+650 pp.

\bibitem{cit_15000_Sawyer}
Sawyer, W. W. Differential equations with polynomial solutions. Quart. J. Math., Oxford Ser. 20, (1949). 22--30.


\bibitem{cit_20000_Suetin}
Suetin, P. K. Classical orthogonal polynomials.
Third edition. Fizmatlit, Moscow, 2005. 480 pp. (Russian)


\bibitem{cit_50000_Gabor_Szego}
Szeg\"o, G\'abor. Orthogonal polynomials. Fourth edition. 
American Mathematical Society, Colloquium Publications, Vol. XXIII. American Mathematical Society, Providence, R.I., 1975. xiii+432 pp.

%\bibitem{cit_70000_T}
%Teschl, Gerald. Jacobi operators and completely integrable nonlinear lattices. Mathematical Surveys and Monographs, 72. 
%American Mathematical Society, Providence, RI, 2000. xvii+351 pp.

\bibitem{cit_80000_Wilkinson}
Wilkinson, J. H. The algebraic eigenvalue problem Clarendon Press, Oxford 1965 {\rm xviii}+662 pp.

\bibitem{cit_90000_Z}
Zagorodnyuk, Sergey M. On generalized Jacobi matrices and orthogonal polynomials. New York J. Math. 9 (2003), 117--136 (electronic).

%\bibitem{cit_92000_Z_banded_matrices}
%Zagorodnyuk, S. M. The direct and inverse spectral problems for some banded matrices.
%Serdica Math. J., {\bf 37} (2011), 9--24.

\bibitem{cit_94000_Z}
Zagorodnyuk S.M., Orthogonal polynomials on rays: properties of zeros,
related moment problems and symmetries.---{\it Journal of Mathematical Physics, Analysis, Geometry}, {\bf 4},
No.~3 (2008), 395-419.

\bibitem{cit_95000_Z}
Zagorodnyuk, S.M.
Orthogonal polynomials related to some Jacobi-type pencils. Ukrain. Mat. Zh.  {\bf 68}, no.~9 (2016), 1180--1190.


\bibitem{cit_97000_Z}
Zagorodnyuk S.M., The inverse spectral problem for Jacobi-type pencils.---{\it SIGMA Symmetry Integrability Geom. Methods Appl.}, {\bf 13}, 
Paper No. 085 (2017), 16 pp.

\bibitem{cit_98000_Zhedanov_JAT}
Zhedanov A., Biorthogonal rational functions and the generalized eigenvalue problem.---{\it 
J. Approx. Theory}, {\bf 101}, 
(1999), no. 2, 303--329.



\end{thebibliography}
\end{document}